\documentclass[a4paper,11pt]{amsart}
\usepackage[applemac]{inputenc}
\usepackage{amsthm,amssymb,bbm,a4wide}
\author{Ivan Marin}
\date{December 5th, 2006}
\address{Institut de Mathématiques de Jussieu, 175 rue du Chevaleret, F-75013 Paris}
\email{marin@math.jussieu.fr}
\urladdr{http://www.math.jussieu.fr/\~\ \!\!\!\!\!\! marin}
\title[Strongly indecomposable finite groups]{Strongly indecomposable finite groups}
\newtheorem{theo}[equation]{Theorem}
\newtheorem{lemma}[equation]{Lemma}
\newtheorem{definition}[equation]{Definition}
\newtheorem{prop}[equation]{Proposition}
\newtheorem{cor}[equation]{Corollary}
\numberwithin{equation}{section}

\def\Hom{\mathrm{Hom}}
\def\Ker{\mathrm{Ker}}
\def\Z{\mathbbm{Z}}

\def\cyc#1{\Z/{#1}\Z}

\def\Fi{\mathrm{Fit\,}}
\begin{document}

\begin{abstract} Motivated by examples in infinite group theory,
we classify the finite groups
whose subgroups can never be decomposed as direct products.
\end{abstract}

\maketitle

\section{Introduction}

It has become clear in recent years that many interesting infinite groups
of geometric origin cannot be decomposed as direct products of
simpler nontrivial groups. Such groups are called \emph{indecomposable}. For instance, under suitable hypothesis,
many fundamental groups cannot be decomposed as direct products unless the space itself admits
such a decomposition. Further, for many of these groups,
it has been noticed that the following property holds.

\begin{definition} A group $G$ is called \emph{strongly indecomposable}
if all its finite-index subgroups are indecomposable.
\end{definition}

Among the groups satisfying this property, we find
\begin{itemize}
\item Infinite Coxeter groups (see \cite{PARIS}) 
\item Zariski-dense subgroups of infinite simple connected algebraic groups (see \cite{HARPE})
\item Mapping class groups with trivial center (see \cite{LONG}) 
\end{itemize}
Also note that, by a theorem of Gromov, every finitely-generated
group of polynomial growth is virtually nilpotent, and that
indecomposability of nilpotent groups is detected by their center
(see \cite{HARPE}). This emphasizes the relevance of this notion.

Coxeter groups, as well as variations of spherical-type Artin groups
(like their commutator subgroups, or their quotient modulo center),
admit large centralizers, but can be thought of as Zariski-dense
subgroups of simple connected algebraic groups (see \cite{MARIN,HARPE}). It is easy to show that
Zariski-dense subgroups of such groups over an infinite field
are strongly indecomposable.

Maybe the most enlightening examples however are generalizations
of free groups, such as torsion-free hyperbolic groups. These groups,
which have ``small centralizers'',
fit into the axiomatic class of non-abelian CSA groups. 
We survey this example in section 2.

It is thus a natural question to ask which
\emph{finite} groups have this
property, namely which finite groups admit only indecomposable
subgroups. We did not find the answer in the litterature on finite
groups, so it is the purpose of this note to classify such groups,
building on classical work in this area.

Recall that the generalized quaternion group $Q_n$ for $n \geq 3$ has
order $2^n$ and is defined by
$$
Q_n = < x,y\ | \ x^{2^{n-1}} = 1, y^2 = x^{2^{n-2}}, y^{-1} x y = x^{-1}>
$$
It is readily checked that all subgroups of $Q_n$ are indecomposable, and
it is known that its center $Z(Q_n)$ has order 2 -- contrasting with the
infinite examples given above.

Our main result is then the following.

\begin{theo} \label{maintheo} Let $G$ be a finite group. All subgroups of $G$
are indecomposable if and only if $G$ is of one of the following types :
\begin{enumerate}
\item $G$ isomorphic to $\Z/p^n\Z$ for some prime $p$.
\item $G$ is generalized quaternion, isomorphic to $Q_n$ for $n \geq 3$.
\item $G = \Z/p^{\alpha} \Z \rtimes \Z/q^{\beta} \Z$
with $p,q$ two different primes, $p$ odd, such that $q^{\beta}$ divides $p-1$,and the image of $\Z/q^{\beta} \Z$
in $(\Z/p^{\alpha} \Z)^{\times}$ has order $q^{\beta}$.
\end{enumerate}
\end{theo}

\section{Examples of strongly indecomposable groups : non-abelian CSA groups}

This class is defined by
the following conditions.

\begin{definition} A group $G$ is said to be CSA if,
for any maximal abelian subgroup $A$ of $G$ and all
$g \in G \setminus A$ we have $A \cap gAg^{-1} = \{ e \}$
(`` A is malnormal '')
\end{definition}

An equivalent definition is to say that every non-trivial element
of $G$ has normal selfnormalizing centralizer. In particular, CSA groups
are commutative-transitive, meaning that nontrivial elements have
abelian centralizers. It is easily checked that every subgroup
of a CSA group is CSA. We refer to \cite{MYAS,JAL} for other properties
of this class of groups.

In fact the relevant class here is the class of \emph{non-abelian}
CSA groups -- note that every abelian group is CSA by definition.
Non-abelian CSA groups
obviously have trivial center. The following
fact is also an easy consequence of the definition : 

\begin{prop} \label{normalCSA} If $G$ is a non-abelian CSA group, then
its non-trivial normal
subgroups are also non-abelian CSA.
\end{prop}
\begin{proof}
We have to show that, if $H<G$ is normal then
it cannot be abelian. Suppose that $H$ is normal, and that
it is abelian. Then it is contained in an abelian maximal subgroup
$\tilde{H}$, with $\tilde{H} \neq G$ since $G$ is non-abelian. Letting $x \in G \setminus \tilde{H}$ it
follows that
$$\{e \} = \ ^x \tilde{H} \cap \tilde{H} \supset \ ^x H \cap H = H \neq \{ e \},$$  
a contradiction. 
\end{proof}
\begin{cor} A non-abelian CSA group is infinite.
\end{cor}
\begin{proof}
It is a classical fact that a non-abelian CSA group $G$ does not contain
any element of order 2 (see e.g. \cite{MYAS} Remark 7). If $G$ were finite, then it would have odd order.
By the Feit-Thompson
theorem it is thus solvable.
It follows that $G$ admits a nontrivial subnormal abelian subgroup,
which is not possible by proposition \ref{normalCSA}.
\end{proof}

A consequence is that finite-index subgroups of non-abelian CSA groups
are non-abelian CSA.

\begin{prop} \label{findexCSA} Let $G$ be a non-abelian CSA group. Then
its finite-index subgroups are non-abelian CSA.
\end{prop}
\begin{proof} Suppose that there exists an abelian finite-index subgroup
$H$ of $G$. We can assume that $H$ is maximal. Let $g_1 =  e $ and
$$ G = g_1 H \sqcup g_2 H \sqcup \dots \sqcup g_n H
$$
the corresponding partition of $G$ in cosets. For all $g \in G$ there
exists $\sigma_g \in \mathfrak{S}_n$ such that $g g_i H = g_{\sigma_g(i)}H$,
and $g \mapsto \sigma_g$ defines a morphism $G \to \mathfrak{S}_n$. Since
$G$ is infinite, there exists $g \neq e$ such that $\sigma_g = e$.
In particuler $\sigma_g(1) = 1$ hence $^g H = H$, and $g \in H$. Since $G$
is non-abelian, $n \geq 2$ and $\sigma_g(2)=2$. It follows that $gg_2 H = g_2 H$
hence $g_2^{-1} g g_2 \in H \cap H^{g_2}$, contradicting once again
the malnormality of $H$. 
\end{proof}

\begin{cor}
Nonabelian CSA groups are strongly indecomposable.
\end{cor}
\begin{proof}
Let $G$ be non-abelian CSA. By proposition \ref{findexCSA} it is sufficient to
show that, if $A$ and
$B$ are nontrivial normal subgroup of $G$, then
$A \cap B \neq \{ e \}$.
Since $A \cap B \supset (A,B)$
one only needs to show $(A,B) \neq \{ e \}$.
Otherwise
the centralizer in $G$ of any $b \in B \setminus \{ e \}$ would contain $A$.
But this centralizer is abelian since $G$ is CSA, and $A$ is not,
a contradiction. 
\end{proof}

\section{Proof of the main theorem}

Our goal here is to classify finite groups whose subgroups
are all indecomposable. We first make some remarks concerning
abelian subgroups.

An indecomposable abelian finite group has to be elementary abelian,
i.e. isomorphic to $\Z/p^n\Z$ for some prime $p$ and positive integer $n$.
Moreover, if a finite strongly indecomposable group $G$ is
not a $p$-group, then its center has to be trivial. Indeed, $Z(G)$ is
abelian hence has order $p^n$ for some
prime $p$. If $n \geq 1$, since $G$ is not a $p$-group it would have
a nontrivial element $x$ of order coprime to $p$ and $<x,Z(G)>$ would
be a decomposable abelian subgroup of $G$, by the Chinese Remainder Theorem.
Hence $Z(G) = \{ 1 \}$ unless $G$ is a $p$-group.

Now assume that $G$ is a $p$-group. We show that, if $G$ is not
abelian, then it is generalized quaternion. Since $G$ is a $p$-group,
its center is nontrivial and contains an element $x$ of order $p$.
Let assume that there exists $y \in G \setminus <x>$ of order $p$.
Then $<y,x>$ would be an abelian
subgroup of $G$ with two distinct subgroups of order $p$,
which is a contradiction
since such a subgroup has to be cyclic.
It follows that $G$ is a $p$-group with exactly
one subgroup of order $p$. By a well-known characterization
(see \cite{ROB} 5.3.6) it follows that $G$ is either cyclic or generalized quaternion.

In particular, the Sylow $p$-subgroups with $p$ odd are
cyclic. In order to tackle the general case, we thus have to distinguish between
two situations : either the 2-Sylow subgroup of $G$ is cyclic, or
it is generalized quaternion.

Before proceeding to the separate study of these two cases, we
first make a remark.

\begin{lemma} \label{nonnormal} If $G$ is strongly indecomposable and is not a $p$-group
then none of its 2-subgroups of $G$
are normal.
\end{lemma}
\begin{proof}
Indeed, let $N$ be a 2-subgroup of $G$. Such
a subgroup $N$ is included in some 2-Sylow subgroup, which is
either a cyclic 2-group or generalized quaternion,
hence contains only one element $z$
of order 2. Since $N$ contains an element of order 2 it contains $z$.
Since $G$ is
not a 2-group, then there exists $g \in G$ of odd order.
But $N \vartriangleleft G$, hence
the element $g z g^{-1}  \in N$ is the only element of order 2 in $G$, and
$gz=zg$. It follows that
$<g,z>$ is a noncyclic abelian subgroup of $G$, which is a contradiction
because such a subgroup is decomposable.
\end{proof}

It follows that the quaternionic case only concerns non-solvable groups :

\begin{lemma} \label{solvable} If $G$ is solvable and strongly indecomposable,
then its 2-Sylow subgroups are cyclic, unless $G$ is a 2-group.
\end{lemma}
\begin{proof}
We assume that $G$ is not a $p$-group. Let $\Fi G$ be the Fitting subgroup
of $G$. Since $G$ is solvable, $\Fi G \neq \{ e \}$. Moreover $\Fi G$ is a nilpotent subgroup of $G$,
hence a direct product of $p$-groups. But $\Fi G$ is indecomposable, hence it is a $p$-group.
Since $\Fi G \vartriangleleft G$ one has $2 \neq p$ by lemma \ref{nonnormal}. In particular,
$\Fi G$ is cyclic. Let $P$ be a 2-Sylow subgroup of $G$. Then $P$ acts on $\Fi G$ by conjugation.
This action is faithful : if $g \in P$ and $x \in \Fi G \setminus \{ 1 \}$,
then $gxg^{-1} = x$ would imply that $<g,x>$ is an abelian noncyclic subgroup of $G$,
a contradiction. It follows that $P$ embeds into $Aut(\Fi G)$, which is abelian
since $\Fi G$ is cyclic. Then $P$ is abelian, hence cyclic.
\end{proof}

\subsection{The cyclic case}

We will use the Hölder-Burnside-Zassenhaus theorem, abreviated HBZ in the sequel,
as stated in \cite{ROB} 10.1.10, p. 281 :

\begin{theo} (Hölder-Burnside-Zassenhaus) If $G$ is a finite group all
of whose Sylow subgroups are cyclic, then $G$ has a presentation
$$
(*) \ \ \ \ \ \ \  G = <a,b \ | \ a^m=1=b^n, b^{-1} a b = a^r >
$$
where $r^n \equiv 1 \ (\mbox{mod } m)$, $m$ is odd, $0 \leq r < m$, and $m$
and $n(r-1)$ are coprime. Conversely, in a group with such a presentation
all Sylow subgroups are cyclic.
\end{theo}

In the special case where $m = p^{\alpha}$, $n = q^{\beta}$, with $p$
and $q$ primes and $\alpha,\beta \geq 1$, these conditions
mean that $p$ is odd, $p \neq q$ and $G \simeq \Z/p^{\alpha} \Z \rtimes \Z/q^{\beta} \Z$.
Indeed, if $a$ is chosen as a generator of $A = \cyc{p^{\alpha}}$, and
$b$ is a generators of $\cyc{q^{\beta}}$, then the action of $B$ on $A$ is given by
$ ^b a = a^r$ for some $r$. Since $b^n = 1$ we have $r^n \equiv 1$
modulo $m$. Since $p \neq q$ we have $\Hom(\cyc{q^{\beta}},
\cyc{p^{\alpha -1}}) = 0$ hence $r$ does not  belong to the subgroup
of $Aut(A) = (\cyc{p^{\alpha}})^{\times}$ isomorphic
to $\cyc{p^{\alpha-1}}$. This subgroup is the set of all $s \in \cyc{p^{\alpha}}$
such that $s \equiv 1$ modulo $p$. It follows that $r-1$ and $p$
are coprime. As a matter of fact, since $Aut(A) \simeq
\cyc{(p-1)} \times \cyc{p^{\alpha-1}}$, the same argument shows
that $r^{p-1} \equiv 1$ modulo $p^{\alpha}$. 

\begin{prop} \label{condmetacyc} Assume $G = \cyc{p^{\alpha}} \rtimes \cyc{q^{\beta}}$
with $\alpha,\beta \geq 1$, $p$ odd and $p \neq q$. Choose a presentation
of $G$ of the form (*). Then $G$ admits no decomposable subgroup iff
$q^{\beta}$ divides $p-1$ and $r$ has order $q^{\beta}$
in $(\cyc{p})^{\times}$.
\end{prop}

\begin{proof}
Let $A = \cyc{p^{\alpha}}$, $B = \cyc{q^{\beta}}$ and $\Phi : B
\to Aut(A) \simeq \cyc{(p-1)} \times \cyc{p^{\alpha-1}}$ the
morphism defining the semidirect product. Choose generators $a_0 \in A$,
$b_0 \in B$ and define $0 < r < p^{\alpha}$ by
$\Phi(b_0)(a_0) = a_0^{b_0} = a_0^r$. The order
of $\Phi(b_0)$ divides $(p-1)p^{\alpha-1}$ and $q^{\beta}$,
hence $p-1$ and $q^{\beta}$. It is the same as the order of $r$
in $(\cyc{p^{\alpha}})^{\times}$ or $(\cyc{p})^{\times}$.

If $G$ admits no decomposable subgroup then $\Phi$ is injective, otherwise
$(\Ker\, \Phi) A$ would provide a decomposable subgroup. It follows
that $q^{\beta}$ divides $p-1$, and $r$ has order $q^{\beta}$
in $(\cyc{p^{\alpha}})^{\times}$, hence in $(\cyc{p})^{\times}$.

Conversely, let $H < G$ be of the form $A' B'$ with $A' = <a_0^u>$
and $B' = <b_0^v >$. Then $b_0^{-v} a_0^u b_0^v = a_0^{ur^v} = a_0^u$
if and only if $u(r^v - 1) \equiv 0$ modulo $p^{\alpha}$.
Since $0 < u < p^{\alpha-1}$ this implies $r^v \equiv 1$ modulo
$p$, hence $q^{\beta}$ divides $v$. Since $0 \leq v < q^{\beta}$
it follows that $v =0$ and $B' = \{ e \}$.

Let us now take an arbitrary subgroup $H$ of $G$. If it is a $p$-group
or a $q$-group it is a subset of some conjugate of $A$ or $B$,
hence it is indecomposable. Otherwise, let $A'$ be a $p$-Sylow of $H$. It
is a subgroup of the unique $p$-Sylow $A$ of $G$, hence of the form
$<a_0^u>$. Let $\varphi : G \to \cyc{q^{\beta}}$ be the natural projection,
and $\overline{D} = \varphi(H)$. It is clear that $A' = \Ker \varphi_{|H}$,
and $\overline{D}$ is a $q$-group, hence $A'$ admits a complement
$D$, which is a $q$-group, by the Zassenhaus theorem. There exists
$g \in G$ such that $ ^{g^{-1}} D = B' \subset B$ by the Sylow theorems.
It follows that $H = \ ^g(A'B')$ because $A' \vartriangleleft G$
and the result follows from the above discussion.

\end{proof}

\subsection{The quaternionic case}

We now assume that $G$ admits a generalized quaternion
group as 2-Sylow subgroup. If $G$ is a $2$-group we are done,
otherwise we have to show that such a group cannot be strongly
indecomposable. For this we show on induction on the order
of $G$ that, if it is strongly indecomposable, then it is solvable.
This will prove theorem \ref{maintheo} by lemma \ref{solvable}.

We may assume that its 2-Sylow subgoups are generalized quaternion,
otherwise it follows from the HBZ theorem that $G$ is solvable.
Then, by a result of Brauer and Suzuki \cite{BS} it cannot be simple
because no simple group admit generalized quaternion group as 2-Sylow
subgroup. Let $N$ be a maximal normal subgroup of $G$.
By the induction hypothesis $N$ is solvable. Since $G/N$
is simple, one only needs to show that it is abelian.

Assume it is not the case. By the Feit-Thompson theorem then $G/N$
has a nontrivial
2-Sylow subgroup $P$.
Moreover, $G/N$ being non-abelian simple cannot be a
$p$-group and $P \neq G/N$. We denote
by $\pi : G \to G/N$ the natural projection and $H = \pi^{-1}(P) \neq G$.
By the induction hypothesis, since $H \neq G$, $H$ is solvable.
By lemma \ref{solvable} we know that $H$ is one of the groups already
classified.  
One of its quotients
contains a nontrivial 2-Sylow hence $H$ is either a 2-group or it
is metacyclic. If it were a 2-group, then $N \subset H$ would be a normal
nontrivial 2-subgroup of $G$, which has been ruled out by lemma \ref{nonnormal}
since $G$ is not
a $p$-group. It follows that $H$ is metacyclic with presentation (*),
and $n$ is
a power of $2$. Since $a \in H$ has odd order, than $\pi(a) = e$ and
$\pi(H) = <\pi(b)>$ is cyclic. But $P = <\pi(H)>$ cannot be cyclic,
because of the classical fact (see \cite{ROB} 10.1.9, p. 280-281) that the
2-Sylow subgroups of a simple non-cyclic group are not cyclic.
It follows by contradiction that $G/N$ is abelian, hence
$G$ is solvable.

This proves by induction that such groups have to be solvable,
but this contredicts lemma \ref{solvable}
and concludes the proof of theorem \ref{maintheo}.

\medskip

\noindent {\bf Acknowledgements.} I thank M. Cabanes for
a useful hint concerning metacyclic groups, and E. Jaligot
for introducing me to the class of CSA groups.

\end{document}